\documentclass{amsart} 
\usepackage{epsf} 
\usepackage{epstopdf} 
\usepackage{amssymb,latexsym, amsmath, amscd, array, graphicx} 

\usepackage{pdfsync} 
\def\hepsffile{\leavevmode\epsffile} 

\swapnumbers \numberwithin{equation}{section} 

\theoremstyle{plain} 
\newtheorem{mainthm}{Theorem} 
\newtheorem{thm}{Theorem}[section]

\theoremstyle{definition}

\newtheorem{remark}[thm]{Remark}

\def\ch{\protect\operatorname{ch}}

\def\sign{\protect\operatorname{sign}}

\def\term{\protect\operatorname{terms}}

\def\Z{{\mathbb Z}} 
 
\def\R{{\mathbb R}}

\def\1{\hbox{\rm\rlap {1}\hskip.03in{\rom I}}} 
\def\Bbbone{{\rm1\mathchoice{\kern-0.25em}{\kern-0.25em} {\kern-0.2em}{\kern-0.2em}I}} 

\def\wh{\widehat}

\def\ov{\overline} \long
\def\forget#1\forgotten{}

\begin{document} 

\title[Intersections of Loops and the Andersen-Mattes-Reshetikhin Algebra] {Intersections of loops and the Andersen-Mattes-Reshetikhin algebra} 
\author[Cahn \& Chernov]{Patricia Cahn and Vladimir Chernov} \address{Patricia Cahn, Department of Mathematics, University of Pennsylvania, David Rittenhouse Lab. 209 South 33rd Street, Philadelphia, PA 19104-6395 USA} \email{pcahn@math.upenn.edu} \address{Vladimir Chernov, Department of Mathematics, 6188 Kemeny Hall, Dartmouth College, Hanover NH 03755-3551, USA} \email{Vladimir.Chernov@dartmouth.edu} 

\begin{abstract}
	Given two free homotopy classes $\alpha_1, \alpha_2$ of loops on an oriented surface, it is natural to ask how to compute the minimum number of intersection points $\#(\alpha_1, \alpha_2)$ of loops in these two classes. 
	
	We show that for $\alpha_1\neq\alpha_2$ the number of terms in the Andersen-Mattes-Reshetikhin Poisson bracket of $\alpha_1$ and $\alpha_2$ is equal to $\#(\alpha_1, \alpha_2)$. Chas found examples showing that a similar statement does not, in general, hold for the Goldman Lie bracket of $\alpha_1$ and $\alpha_2$. 
	
	The main result of this paper in the case where $\alpha_1, \alpha_2$ do {\bf not} contain different powers of the same loop first appeared in the unpublished preprint of the second author. In order to prove the main result for all pairs of $\alpha_1\neq \alpha_2$ we had to use the techniques developed by the first author in her study of operations generalizing Turaev's cobracket of loops on a surface.
\end{abstract}

\maketitle
\section{Basic definitions and main results} 

We work in the smooth category, where ``smooth'' means $C^{\infty}$. Throughout this work $F$ is an oriented two dimensional surface and $\widehat \pi$ is the set of free homotopy classes of loops on $F$. For $\beta\in \pi_1(F)$ we let $\langle \beta \rangle \in \widehat \pi$ denote the free homotopy class corresponding to $\beta$. Similarly for a loop $a:S^1\to F$ we denote by $\langle a\rangle \in \widehat \pi$ the free homotopy class containing $a$. We let $\Z[\widehat \pi]$ be the free abelian group of formal finite linear combinations of the elements of $\widehat \pi$ with integer coefficients.

Goldman~\cite{Goldman} constructed a Lie algebra structure on $\Z[\widehat \pi]$.  Turaev~\cite{Turaevskein} later defined a cobracket on the free $\Z$-module generated by the set of {\it nontrivial} free homotopy classes, and showed that his cobracket, together with the Goldman bracket, forms a Lie bialgebra structure on this module.  Chas and Sullivan's String Topology generalizes the Goldman-Turaev Lie bialgebra \cite{ChasSullivan1, ChasSullivan2}.   

Given an element $\gamma\in \Z[\widehat \pi]$ with $\gamma=\sum_k j_k\omega_i$, $j_k\in \Z, \omega_k\in \widehat \pi$ we put $\term(\gamma)=\sum_k|j_k|$ and call it the {\it number of terms of $\gamma$.\/}  We always assume that the expression $\sum_k j_k\omega_i$ is reduced.

Given two smooth loops $a_1, a_2:S^1\to F$ we let $i(a_1, a_2)$ be the (possibly infinite) number of intersection points of $a_1$ and $a_2$, i.e., the number of pairs $(t_1, t_2)\in S^1\times S^1$ such that $a_1(t_1)=a_2(t_2).$ For $\alpha_1, \alpha_2\in \widehat \pi$ we put 
\begin{equation}
	\#(\alpha_1, \alpha_2)=\min_{a_1\in \alpha_1, a_2\in \alpha_2} i(a_1, a_2) 
\end{equation}
and we call it the {\it minimal intersection number} of the free homotopy classes $\alpha_1$ and $ \alpha_2.$

It is easy to see that $\term([\alpha_1, \alpha_2])\leq \#(\alpha_1, \alpha_2),$ where $[\cdot, \cdot]$ denotes the Goldman Lie bracket of two elements of $\Z[\widehat \pi].$ Goldman~\cite{Goldman} showed that if $\alpha_1$ contains a simple loop and $[\alpha_1, \alpha_2]=0$, then $\#(\alpha_1, \alpha_2)=0$ and thus there exist $a_1\in \alpha_1, a_2\in \alpha_2$ such that $a_1$ and $a_2$ do not intersect. 

Later Chas~\cite[Corollary, page 27]{ChasMinimalIntersection} proved that the following two statements are equivalent:
\begin{description}
	\item[I] $\alpha_1$ contains a power of simple loop, 
	\item[II] $\term([\alpha_1, \alpha_2])=\#(\alpha_1, \alpha_2)$ for all $\alpha_2\in \widehat \pi$. 
\end{description}

Thus if $\alpha_1$ contains a power of a simple loop, then for every $\alpha_2$ one can use Goldman's Lie bracket to compute the minimal number of intersection points of loops in homotopy classes $\alpha_1$ and $\alpha_2$. On the other hand, if $\alpha_1$ does not contain a power of a simple loop, then for some $\alpha_2$ we have $\term([\alpha_1, \alpha_2])<\#(\alpha_1, \alpha_2)$. In fact, it is even possible to have $\term([\alpha_1, \alpha_2])=0$ while $\#(\alpha_1, \alpha_2)\neq 0$. The first such examples were found by Chas~\cite[Example 5.6]{ChasCombinatorial}.

Andersen, Mattes and Reshetikhin~\cite{AMR1},~\cite{AMR2} constructed a Poisson algebra on a $\mathbb{Z}-$module generated by chord diagrams on $F$. An element $\alpha\in \widehat \pi$ can be viewed as a chord diagram with one circle and zero chords. For $\alpha_1, \alpha_2\in \widehat \pi$ their Andersen-Mattes-Reshetikhin Poisson bracket $\{\alpha_1, \alpha_2\}=\sum_k j_k\tau_k$ is a finite integer linear combination of chord diagrams $\tau_k$ consisting of two circles connected by a chord. We put $\term(\{\alpha_1, \alpha_2\})=\sum_k|j_k|$ and call it {\it the number of terms of $\{\alpha_1, \alpha_2\}$.\/} We always assume that the expression $\sum_k j_k\tau_k$ is reduced.

The main result of this work is the following theorem.
\begin{mainthm}
	\label{maintheorem} Let $\alpha_1\neq \alpha_2$ be two free homotopy classes of loops on an oriented surface $F.$ Then $\term(\{\alpha_1, \alpha_2\})=\#(\alpha_1, \alpha_2)$, i.e.~the number of terms in the Andersen-Mattes-Reshetikhin Poisson bracket of $\alpha_1$ and $\alpha_2$ equals the minimum number of intersection points of loops in the classes $\alpha_1$ and $\alpha_2.$ 
\end{mainthm}

The proof of this theorem is contained in Section~\ref{proofofmaintheorem}.

Since $\{\cdot, \cdot\}$ is skew symmetric and hence $\{\alpha, \alpha\}$ is always zero, the statement of this theorem is false for the case when $\alpha_1=\alpha_2.$ However the conditions of the theorem allow the case where $\alpha_1=\langle \delta^i\rangle $ and $\alpha_2=\langle \delta^j\rangle$ for $\delta\in \pi_1(F)$ and $i\neq j$, i.e.~when $\alpha_1$ and $\alpha_2$ contain different powers of the same loop. 

Of course when $\alpha_1=\alpha_2$ we have that $\#(\alpha_1, \alpha_2)=\#(\alpha_1, \alpha_1^{-1}),$ where $\alpha_1^{-1}$ denotes the class containing the loops from $\alpha_1$ with orientation reversed. By Theorem~\ref{maintheorem} we have $\#(\alpha_1, \alpha_1^{-1})=\term(\{\alpha_1, \alpha_1^{-1}\})$ and we see that $\#(\alpha_1, \alpha_2)$ can easily be computed using the Andersen-Mattes-Reshetikhin Poisson bracket even in the case where $\alpha_1=\alpha_2.$

Given a smooth loop $a:S^1\to F$ we let $i(a)$ be the (possibly infinite) number of self-intersection points of $a$, i.e., one half of the number of pairs $(t_1, t_2)\in S^1\times S^1$ with $t_1\neq t_2$ such that $a(t_1)=a(t_2).$ 
For $\alpha \in \widehat \pi$ we put 
\begin{equation}
	\#(\alpha)=\min_{a\in \alpha} i(a) 
\end{equation}
and we call it the {\it minimal self-intersection number} of the free homotopy class $\alpha.$

There is a simple relationship between the minimum number of self-intersection points $\#(\alpha)$ of a generic loop in $\alpha$ and the minimal intersection number $\#(\alpha,\alpha)=\#(\alpha,\alpha^{-1})$; we will see that $\#(\alpha,\alpha)=2(\#(\alpha)-(n-1))$ where $n$ is the largest integer such that $\alpha=\beta^n$ for some class $\beta\in \pi_1(F)$.  This gives rise to the following formula for the minimal self-intersection number $\#(\alpha)$.

\begin{mainthm} \label{mainthmcor}
	Let $\alpha$ be a nontrivial free homotopy class of loops on an oriented surface $F$, and let $n$ be the largest integer such that $\alpha=\beta^n$ for some class $\beta\in \pi_1(F)$.  Let $p$ and $q$ be distinct nonzero integers.  Then $\#(\alpha)=\frac{1}{2|pq|}\term(\{\alpha^p,\alpha^{q}\})+n-1$.
\end{mainthm}

Similar formulas for the self-intersection number in terms of the Goldman bracket  of $\alpha^p$ and $\alpha^q$ were given by Chas and Krongold \cite{ChasKrongold}, in the case where the surface has nonempty boundary, $n=1$, $p$ and $q$ are positive, and $p$ or $q$ is at least 3.  It also is interesting to note the resemblance between the formula in Theorem \ref{mainthmcor} and the formula $\#(\alpha)=\frac{1}{2}\term \mu(\alpha)+n-1$ given by the first author in \cite{CahnSelfIntersection}, where $\mu$ is an operation generalizing Turaev's cobracket.  
\begin{remark}
	{\it Garlands and some previously known cases of Theorem~\ref{maintheorem}.\/} In the unpublished preprint of the second author~\cite[Theorem 1.4]{ChernovPreprintIntersection} we introduced a graded Poisson algebra on the oriented bordism group of mappings into a manifold $M$ of garlands made out of a collection $\mathfrak N$ of odd dimensional manifolds. In the case when garlands are made out of circles and $M$ is a two dimensional surface, zero-dimensional bordisms form a Poisson subalgebra. We showed~\cite[Proposition 12.4 and Section 12.2]{ChernovPreprintIntersection} that this subalgebra is identified with a symmetrization of the subalgebra of the Andersen-Mattes-Reshetkhin Poisson algebra formed by tree-like chord diagrams. The study of garland algebras was initiated in the joint work of Rudyak and the second author~\cite{ChernovRudyakGarlands}.
	
	In~\cite[Theorem 1.8, Theorem 10.1 and Section 10.2]{ChernovPreprintIntersection} we provided the sketch of the proof of the statement of Theorem~\ref{maintheorem} under the assumption that it is {\bf not} true that $\alpha_1=\langle \delta^i\rangle $ and $\alpha_2=\langle \delta^j\rangle$ for some $\delta\in \pi_1(F)$ and $i,j\in\Z.$ The technique used in~\cite{ChernovPreprintIntersection} does not seem to allow us to prove Theorem~\ref{maintheorem} in the case where $\alpha_1$ and $\alpha_2$ contain different powers of the same element of $\pi_1(F).$
	
	The proof of Theorem~\ref{maintheorem} in this paper uses the techniques developed in the work of the first author~\cite{CahnSelfIntersection}. In that work it is proved~\cite[Theorem 1.1 and Theorem 1.2]{CahnSelfIntersection} that if one generalizes the Turaev cobracket in the spirit of the Andersen-Mattes-Reshetikhin Poisson algebra, then the resulting operation gives a formula for the minimum number of self intersection points of a loop in any given $\alpha\in \widehat \pi.$ 
	
	We also use some of the ideas from the papers of Turaev~\cite{Turaev} and 
of Turaev-Viro~\cite{TuraevViro}, who constructed algebraic operations that compute the minimal number of self-intersection points and the minimal number of intersection points of loops in given free homotopy classes. We are not aware of any relation between Turaev-Viro operations and the Poisson algebra of Andersen-Mattes-Reshetikhin. 
	
	Our results allow one to compute the minimal intersection number, and we give an example in the last section of the paper which illustrates this for surfaces with boundary.  However, we do not place much emphasis on algorithmic methods of computing the minimal intersection number in this paper.  For algorithmic methods of computing the minimal intersection number, we refer the reader to the works of Lustig~\cite{Lustig} and Cohen-Lustig~\cite{CohenLustig}.

\end{remark}

\section{The Goldman Lie Algebra}
Goldman~\cite{Goldman} introduced a Lie algebra on $\Z[\wh \pi].$ 
Given two free loop homotopy classes $\alpha_1, \alpha_2\in \wh \pi,$ their Goldman's Lie bracket 
$[\alpha_1, \alpha_2]\in \Z[\wh \pi]$ is defined as follows. 
Take two transverse representatives $a_1, a_2$ of $\alpha_1, \alpha_2$ that are generic in the sense that 
intersection points of $a_1$ and $a_2$ are transverse double points and they do not happen at the self
intersection points of $a_1$ or $a_2.$ 

Put $P=a_1\cap a_2.$ For each $p\in P$ we construct a two-frame at $T_pF$ whose first vector is the velocity 
vector of $a_1$ at $p$ and whose second vector is the velocity vector of $a_2$ at $p.$
Put $\sign(p)$ to be $+1$ if the orientation of $T_pF$ given by the 
above frame is positive and put $\sign(p)=-1$ if the orientation of $T_pF$ given by the above frame is 
negative. 

After an orientation preserving reparameterization, the loops 
$a_1$ and $a_2$ can be considered as loops based at $p.$ We denote by  $\langle a_1\cdot_p a_2\rangle$
the element of $\wh \pi$ realized by the class of the product $a_1a_2\in \pi_1(F,p)$. 

{\it Goldman Lie bracket $[\alpha_1, \alpha_2]$ is defined by \/} 
\begin{equation}\label{goldman}
[\alpha_1, \alpha_2]=\sum_{p\in P}\sign (p)\langle a_1\cdot_p a_2\rangle\in \Z[\wh \pi],
\end{equation}
and it can be extended to $\Z[\wh \pi]$ by linearity.

Turaev~\cite{Turaevskein} introduced a cobracket on the quotient of $\Z[\wh \pi]$ by the $\Z$-submodule 
generated by the class of the trivial loop. Together Goldman's bracket and the Turaev's cobracket make this 
quotient into an involutive Lie bialgebra.

\section{The Andersen-Mattes-Reshetikhin Poisson Algebra}\label{AMRalgebra} We now summarize the construction of the Andersen-Mattes-Reshetikhin Poisson algebra given in ~\cite{AMR1},~\cite{AMR2}.

{\it A chord diagram\/} is a topological space that consists of some number of disjoint oriented circles $S^1_i, i=1, \cdots, q,$ and disjoint arcs $A_j, j=1, \cdots, r,$ such that the end points of the arcs are distinct and $\cup _j
\partial A_j=(\cup _i S^1_i)\cap (\cup_j A_j).$ The circles $S^1_i$ are called {\it core components\/} and the arcs $A_j$ are called {\it chords\/} of the diagram.

A {\it geometric chord diagram on an oriented surface $F$\/} is a smooth map of a chord diagram to the surface that maps each chord to a point. {\it A generic geometric chord diagram on $F$\/} is a geometric chord diagram such that all the circles are immersed and all the multiple points between them are transverse double points. A {\it chord diagram on $F$\/} is an equivalence class of geometric chord diagrams modulo homotopy. The commutative {\it multiplication of two chord diagrams on $F$\/} is defined to be their union. 

\begin{figure}
	[htbp] 
	\begin{center}

		\epsfxsize 10cm 
		\hepsffile{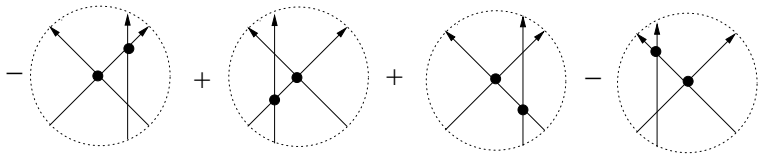} 
	\end{center}
	\caption{$4T$-relation. The diagrams are identical outside of the dotted circles, and the parts that are drawn are embedded in $F.$}\label{4Trelation.fig} 
\end{figure}

Consider the complex vector space $X$ whose basis is the set of chord diagrams on $F.$ Let $Y$ be the subspace of $X$ generated by linear combinations called {\it $4T$-relations}. One of the $4T$-relations is depicted in Figure~\ref{4Trelation.fig}. The others are obtained by reversing the orientations of strands in Figure~\ref{4Trelation.fig}, by following the rule that for each chord that intersects a component whose orientation is reversed we get a factor of $(-1)$ in front of the diagram.

The subspace $Y$ is an ideal with respect to multiplication, and the {\it Andersen-Mattes-Reshetikhin chord algebra} is $\ch(F)=X/Y.$

Given two chord diagrams $[D_1], [D_2]$ on $F,$ pick two geometric chord diagrams $D_1,D_2$ representing them so that $D_1\cup D_2$ is a generic chord diagram. Put $P=D_1\cap D_2$ to be the set of intersection points of $D_1$ and $D_2.$ Each $p\in P$ is the intersection point of a curve from $D_1$ with a curve from $D_2.$ For each $p\in P$ we construct a $2$-frame in $T_pF$ whose first vector is the velocity vector of $D_1$ at $p$ and whose second vector is the velocity vector of $D_2$ at $p$. We put $\sign(p)=+1$ if this $2$-frame gives a positive orientation of $T_pF,$ and we put $\sign(p)=-1$ otherwise. For $p\in P$ we put $D_1\cup_p D_2$ to be the chord diagram on $F$ obtained by joining $D_1^{-1}(p)$ and $D_2^{-1}(p)$ by a chord.

{\it The Andersen-Mattes-Reshetikhin Poisson bracket $\{\cdot, \cdot\}$ of two chord diagrams $[D_1], [D_2]$ on $F$\/} is given by 
\begin{equation}
	\label{AMRPoisson} \{[D_1], [D_2]\}=\sum _{p\in D_1\cap D_1}\sign(p) [D_1\cup_p D_2], 
\end{equation}
where a geometric diagram in square brackets denotes the chord diagram equivalence class realized by it. 

Note that the Andersen-Mattes-Reshetikhin Poisson bracket can be viewed as a generalization of Goldman's bracket.  Namely, if one computes $\{\alpha_1,\alpha_2\}$ for two free homotopy classes, and then smoothes the resulting diagrams at their chords according to the orientations on the core circles, one recovers the Goldman bracket $[\alpha_1,\alpha_2]$ of the two loops.  Thus the bound on the minimal intersection number given by the Andersen-Mattes-Reshetikhin bracket must be at least as good as the bound given by Goldman's bracket.  Our main theorem says that the first bound is in fact an equality, while the results of Chas~\cite{ChasCombinatorial} say that in general the second bound is not an equality. 
\section{Proof of Theorem~\ref{maintheorem}}\label{proofofmaintheorem} 

If $\alpha_1$ contains contractible loops, then $\{ \alpha_1, \alpha_2\}=0$. Thus $\term(\{ \alpha_1, \alpha_2\})=0.$ Furthermore, two loops from the homotopy classes $\alpha_1, \alpha_2$ can be made disjoint by a homotopy. Thus $\#(\alpha_1, \alpha_2)=0=\term(\{ \alpha_1, \alpha_2\})$, so from now on we assume that $\alpha_1$ and $\alpha_2$ do not contain contractible loops.  

{\bf We begin by assuming that $F$ is a compact surface.\/} (Later, we will consider the case where $F$ is not compact.) If $F$ is $S^2$ or the annulus $A$, then every two loops on $F$ can be made disjoint by homotopy. Thus $\#(\alpha_1, \alpha_2)=0$. We also have $\term (\{\alpha_1, \alpha_2\})=\term(0)=0$, which proves the theorem for $F=S^2, A$.

{\bf First we consider the case where $F^2$ is the torus $T^2.$\/} Clearly $\term (\{\alpha_1, \alpha_2\})$ is at least the absolute value of the intersection number $[\alpha_1]\bullet [\alpha_2]$ of the classes in $H_1(F^2)$ realized by $\alpha_1$ and $\alpha_2.$ Since $H_1(T^2)=\pi_1(T^2)=\Z\oplus\Z,$ we get that $\alpha_1=i_1(m_1, l_1)$ for some $i_1\in \Z$ and coprime $m_1, l_1\in \Z.$ Take $m_2, l_2\in \Z$ such that $m_1l_2-l_1m_2=1.$ Then $\alpha_2=i_2(m_1, l_1)+j_2(m_2, l_2),$ for some $i_2, j_2\in \Z.$ Put $\{{\bf i, j}\}$ to be the standard orthonormal frame in $\R^2$ and take $\mu$ to be the simple loop on $T^2$ that lifts to a straight line $\vec {\bf r_1}(t)=m_1t{\bf i}+l_1t {\bf j}$ in the total space $\R^2$ of the universal cover $\R^2\to \R^2/(\Z \oplus\Z)=T^2.$ Take $a_1\in \alpha_1$ to be $\mu^{i_1.}$

Put $\ov \mu$ to be a small parallel shift of $\mu$ and put $\lambda$ to be a simple loop such that $\lambda$ and $\ov \mu$ have the same base point and such that $\lambda$ lifts to a straight line ${\bf r_2}(t)=m_2t {\bf i}+l_2t{\bf j}$ in the universal cover $\R^2$ of $T^2.$ Put $a_2=\ov \mu_1^{i_2}\lambda^{j_2}\in \alpha_2.$ It is easy to see that the number of intersection points of $a_1$ and $a_2$ is $|i_1j_2|$ and that $i_1j_2=[\alpha_1]\bullet [\alpha_2].$

Thus $\#(\alpha_1, \alpha_2)\leq |[\alpha_1]\bullet [\alpha_2]|.$ Clearly $|[\alpha_1]\bullet [\alpha_2]| \leq \#(\alpha_1, \alpha_2)$ and so $\#(\alpha_1, \alpha_2)=|[\alpha_1]\bullet [\alpha_2]|.$

Since $\term(\{\alpha_1, \alpha_2\})\leq \#(\alpha_1, \alpha_2)=|[\alpha_1]\bullet [\alpha_2]|$ and $\term (\{\alpha_1, \alpha_2\})\geq |[\alpha_1]\bullet [\alpha_2]|$ we have that $\#(\alpha_1, \alpha_2)=\term(\{\alpha_1, \alpha_2\})$ and this finishes the proof of Theorem~\ref{maintheorem} for $F=T^2.$

{\bf Next we consider the case where $F$ is a compact surface, with or without boundary, other than $S^2, T^2$, or $A$\/}. After this, we will consider the case where $F$ is a noncompact surface. 

We use a decomposition of $F$ into pairs of pants to construct a Riemannian metric $g$ of constant negative curvature $-1$ on $F$ such that the boundary of $F$ is geodesic. Throughout the proof we will use the following two facts about $\pi_1(F)$:

{\it 1. All nontrivial abelian subgroups of $\pi_1(F)$ for such $F$ are infinite cyclic, and moreover for each $\alpha\neq 1\in \pi_1(F)$ there is a unique maximal infinite cyclic subgroup containing $\alpha.$ } For closed $F$ this follows from the 
Preissman Theorem~\cite[Theorem 3.2 page 260]{docarmo},~\cite{Preissman}. For compact $F$ with 
$\partial F\neq \emptyset$ this holds since $\pi_1(F)$ is free.

{\it 2. Every nontrivial free homotopy class on $(F,g)$ is realizable by a geodesic loop that is unique up to reparameterization;\/} see~\cite[Lemma B.4.5]{BenedettiPetronio} for the case $
\partial F=\emptyset,$ or see~\cite[Theorem 1.6.6]{Buser} for the general case. For closed $F$ this is a result of Hadamard~\cite{Hadamard} which is a particular case of the Cartan Theorem~\cite{Cartan},~\cite[Theorem 3.8.14]{Klingenberg} that holds for closed negative sectional curvature manifolds of all dimensions. 

Let $a_i$ be the closed geodesic realizing $\alpha_i, i=1,2.$ Let $\beta_i, i=1,2,$ be the generator of the maximal infinite cyclic group containing the class of $a_i\in \pi_1(F^2, a_i(1)),$ so that $\alpha_i=\langle \beta_i^{k_i}\rangle $ for some $k_i\in \Z.$ Choosing $\beta_i^{-1}$ instead of $\beta_i,$ if needed, we can assume that $k_i>0.$

The free homotopy class $\langle \beta_i \rangle$ is realizable by a unique (up to reparameterization) closed geodesic $b_i$ $i=1,2.$ Thus $\alpha_i$ is realizable by a closed geodesic that travels along the oriented closed geodesic $b_i$ a total of $k_i$ times. Furthermore this choice of $a_i$ is the unique closed geodesic realizing $\alpha_i$. Observe that $\langle b_i \rangle \in \pi_1(F^2)$ is not realizable as a power of another element. Hence $b_i$ is injective on the complement of a finite set of points in $S^1$ that are the preimages of the transverse self-intersection points of $b_i, i=1,2.$

If one of $b_1$, $b_2$ is a boundary curve and $b_1\neq b_2^{\pm 1}$ then $\{\alpha_1, \alpha_2\}=\{\langle b_1^{k_1}\rangle, \langle b_2^{k_2}\rangle\}=0$, so $\term\{\alpha_1, \alpha_2\}=0.$ On the other hand, $i(b_1^{k_1}, b_2^{k_2})=0$, and hence $\#(\alpha_1, \alpha_2)=0.$ Thus the statement of the theorem holds in this case. If $b_1=b_2^{\pm 1}$ then we can take a representative of $\alpha_1$ located in the thin neighborhood of the boundary component and such that this representative does not touch the boundary. Taking a representative of $\alpha_2$ located in an even thinner neighborhood of the boundary we get that $\#(\alpha_1, \alpha_2)=0=\term(\{\alpha_1, \alpha_2\}).$

Now in the remainder of the proof we may assume that the geodesics $b_1, b_2$ are located in $F\setminus 
\partial F.$

In Subsection~\ref{powersofdifferentelements} we prove the theorem in the cases where $\alpha_1$ and $\alpha_2$ do not contain powers of the same loop. This means that $b_1$ and $b_2$ intersect at a finite number of points and these intersections are transverse. In Subsection~\ref{powersofthesameelement} we prove the theorem in the case where $\alpha_1$ and $\alpha_2$ contain powers of the same loop $\delta$.

To present the proof we will need to introduce the chord diagrams $b_1\bullet_p b_2$ and $b_1^{k_1}\bullet_p b_2^{k_2}$, and we do this in the following subsection.

\subsection{$b_1\bullet_p b_2$ and its properties}\label{bulletanditsproperties} 
Note that in this section and in the remainder of the paper, we abuse notation slightly and refer to the homotopy class of a path or a loop $b$ by $b$ rather than $\langle b \rangle$.  It should be clear from context whether we are referring to a loop or path, or to the class of a loop or the class of a path.  Assume that $b_1, b_2$ are two transverse immersed closed curves that are not necessarily in general position. This means that their intersection points could have multiplicity greater than two. (In our proof, the $b_i$ are closed geodesics, so we cannot assume all intersection points are double points.) For $p', p''\in S^1$ such that $b_1(p')=b_2(p''),$ we let $p=b_1(p')=b_2(p'')$ and we let $\sign(p)$ be the sign of the intersection point of the branches of $b_1$ and $b_2$ containing $p$ and $p'$ respectively. 

Note that this is a slight abuse of notation since $p$ could be an intersection point of multiplicity greater than $2$. However, it will be clear from context which pair of preimages of an intersection point of $b_1$ and $b_2$ we are talking about. 

We let $b_1 \bullet _p b_2$ be the homotopy class of a geometric chord diagram with two circles and one chord connecting them, such that the map on one circle is $b_1$, the map on the other circle is $b_2$, and the chord connects the points $p'$ and $p''$ on the two circles.

For nonzero integers $k_1, k_2$, we let $b_1^{k_1}\bullet _p b_2^{k_2}$ be the homotopy class of a geometric chord diagram with two circles and one chord connecting them, such that the map on one circle is $b_1^{k_1}$, and the map on the second circle is $b_2^{k_2}$. We can view $b_i^{k_i}$ as a composition of $b_i$ with the covering map $S^1\rightarrow S^1,$ $e^{2\pi i t}\to e^{k_i2 \pi i t}$ of degree $k_i$. Recall that $p$ is an intersection point of $b_1$ and $b_2$ with preimages $p'$ and $p''$. The points $p'$ and $p''$ have $k_1$ and $k_2$ preimages each under the covering maps of degree $k_i$. The chord of the diagram $b_1^{k_1}\bullet _p b_2^{k_2}$ connects any pair of preimages of the points $p'$ and $p''$ on the two circles under the covering maps. Note that the homotopy class $b_1^{k_1}\bullet _p b_2^{k_2}$ does not depend on the choices of the preimages of $p'$ and $p''$ that are to be connected by a chord.

Assume that $b_1^{k_1}\bullet _p b_2^{k_2}$ and $b_3^{k_3}\bullet _p b_4^{k_4}$ are homotopy classes of two chord diagrams for which the chord is mapped to the same point $p$ on the surface. It is easy to see that these homotopy classes of chord diagrams are equal if and only if there exists $s\in \pi_1(F,p)$ such that one of the following conditions holds 
\begin{equation}\label{descr}
\begin{split}
	{\bf 1.\/}\, \text{ }s(b_1^{k_1})s^{-1}=b_3^{k_3}\in\pi_1(F,p) \text{ and }s(b_2^{k_2})s^{-1}=
b_4^{k_4}\in\pi_1(F,p)\\ 
	{\bf 2.\/}\, \text{ }s(b_2^{k_2})s^{-1}=b_3^{k_3}\in\pi_1(F,p) \text{ and }s(b_1^{k_1})s^{-1}=
b_4^{k_4}\in\pi_1(F,p). 
\end{split}
\end{equation}

That is these two homotopy classes of chord diagrams are equal exactly when the two loops in one diagram, taken in either order, can be simultaneously conjugated by some element of $\pi_1(F)$ to obtain the two loops in the other diagram.

{\it Now we see we only need Part ${(\bf 1.)}$ of Equation \ref{descr}.}  In Theorem \ref{maintheorem}, we assume $\alpha_1$ and $\alpha_2$ are distinct free homotopy classes.  A term of $\{\alpha_1,\alpha_2\}$ is of the form $a\bullet_p b$ where either $a\in \alpha_1$ and $b\in \alpha_2$ or $a\in\alpha_2$ and $b \in \alpha_1$.  During the proof of Theorem 1 we will always assume that the terms are written so that $a\in \alpha_1$ and $b\in \alpha_2$. Now consider two terms of $\{\alpha_1,\alpha_2\}$, say $a\bullet _p b$ and $c\bullet_p d$, corresponding to the same intersection point $p$.  Clearly Part ${(\bf 2.)}$ of Equation \ref{descr} cannot hold, because $b$ and $c$ are not in the same free homotopy class, and $a$ and $d$ are not in the same free homotopy class.  Therefore if $a\bullet _p b$ and $c\bullet_p d$ are homotopic chord diagrams, they are related by Part ${(\bf 1.)}$ of Equation \ref{descr}.

\subsection{The proof of Theorem~\ref{maintheorem} in the case where $\alpha_1$ and $\alpha_2$ do {\bf not\/} contain powers of the same loop.}\label{powersofdifferentelements}

It is easy to observe that for such $\alpha_1, \alpha_2$ we have $\{\alpha_1, \alpha_2\}=k_1k_2\sum_{p\in b_1\cap b_2}\sign (p) (b_1^{k_1}\bullet _p b_2^{k_2}).$ 

To prove the Theorem it suffices to show that for $p_1\neq p_2$ in $b_1\cap b_2$ we have that the geometric chord diagrams $b_1^{k_1}\bullet _{p_1} b_2^{k_2}$ and $b_1^{k_1}\bullet _{p_2} b_2^{k_2}$ are not homotopic. 

Recall that the curves $b_1, b_2$ could have intersection points of multiplicity greater than two, so here and below by $p_1\neq p_2$ we mean that the two pairs of preimages of the intersection points $p_1$ and $p_2$ that are hidden in the definition of $p_1$ and $p_2$ are different.

We argue by contradiction and assume that the chord diagrams $b_1^{k_1}\bullet _{p_1} b_2^{k_2}$ and $b_1^{k_1}\bullet _{p_2} b_2^{k_2}$ are homotopic.

\begin{figure}
	[htbp] 
	\begin{center}

		\epsfxsize 11cm \hepsffile{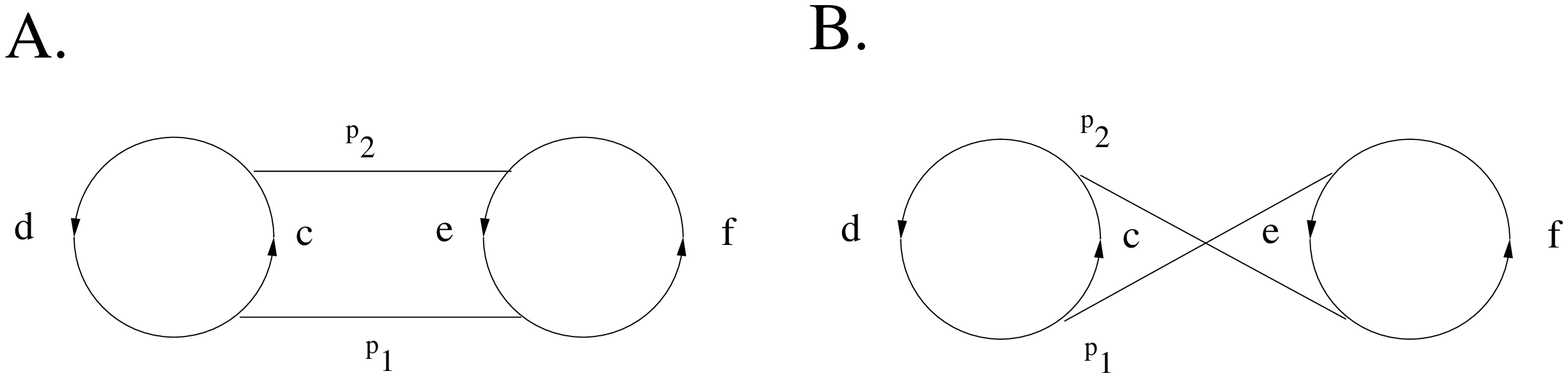} 

	\end{center}
	\caption{Two possible mutual positions for chords corresponding to intersection points of $b_1$ and $b_2$}\label{chorddiagrams.fig} 
\end{figure}

At first it may seem that there are two possible cases for the mutual position of chords corresponding to $p_1$ and $p_2.$ They are shown in ~\ref{chorddiagrams.fig}.A and Figure~\ref{chorddiagrams.fig}.B, respectively.  However we can turn Figure~\ref{chorddiagrams.fig}.B into ~\ref{chorddiagrams.fig}.A by sliding $p_2$ around the circle and then exchanging the labels $p_1$ and $p_2.$ 

 Note that some of the circle arcs between the chord ends on these figures could be trivial. However the two chords are really different, since we assumed $p_1\neq p_2.$

{\bf Here we consider the case described in Figure~\ref{chorddiagrams.fig}.A} In this case we get 
\begin{equation}
	b_1^{k_1}\bullet _{p_1} b_2^{k_2} =(cd)^{k_1}\bullet _{p_1} (fe)^{k_2} 
\end{equation}
and
\begin{equation}
	b_1^{k_1}\bullet _{p_2} b_2^{k_2} = (dc)^{k_1}\bullet _{p_2} (ef)^{k_2} . 
\end{equation}
We first deform the second chord diagram so that its chord is mapped to the point $p_1$: 
\begin{equation}
	(dc)^{k_1}\bullet _{p_2} (ef)^{k_2} = c(dc)^{k_1}c^{-1}\bullet _{p_1} c(ef)^{k_2}c^{-1} = (cd)^{k_1}\bullet _{p_1} c(ef)^{k_2}c^{-1} . 
\end{equation}
Now both chord diagrams can be considered as pairs of loops based at $p_1$; see Equation~\eqref{descr}. Since we assumed that $ b_1^{k_1}\bullet _{p_1} b_2^{k_2} = b_1^{k_1}\bullet _{p_2} b_2^{k_2} $ and $\alpha_1\neq \alpha_2$, we get that there exists $\mu\in \pi_1(F, p_1)$ such that 
\begin{equation}
	\label{case2Aeqn1} \mu(cd)^{k_1}\mu^{-1}=(cd)^{k_1}\in \pi_1(F, p_1) 
\end{equation}
and 
\begin{equation}
	\label{case2Aeqn2} \mu c (ef)^{k_2} c^{-1}\mu^{-1}=(fe)^{k_2}\in \pi_1(F, p_1). 
\end{equation}

Since the centralizer of $(cd)^{k_1}\in \pi_1(F, p_1)$ is the infinite cyclic group generated by $(cd)\in\pi_1(F, p_1)$, we have that $\mu=(cd)^k\in \pi_1(F, p_1)$ for some $k\in\Z.$ Substitute this into equation \eqref{case2Aeqn2} to get 
\begin{equation}
	(cd)^k c (ef)^{k_2} c^{-1}(cd)^{-k}=(fe)^{k_2}=f(ef)^{k_2}f^{-1}\in \pi_1(F, p_1). 
\end{equation}

Thus the paths $f^{-1}(cd)^k c (ef)^{k_2}$ and $(ef)^{k_2}f^{-1}(cd)^kc$ are homotopic as paths with fixed ends at $p_2$. Thus $f^{-1}(cd)^k c$ and $(ef)^{k_2}$ are commuting elements of $\pi_1(F, p_2).$ Since $(ef)\in \pi_1(F,p_2)$ generates the centralizer of $(ef)^{k_2}\in \pi_1(F,p_2)$, we get that $f^{-1}(cd)^k c=(ef)^l\in\pi_1(F, p_2)$ for some $l\in\Z.$ Hence $(cd)^k c$ and $f(ef)^l$ are homotopic as arcs with fixed end points. Now one considers four possible cases of signs of $k$ and $l.$

We consider the case when both $k$ and $l$ are negative. The other cases are considered in a similar fashion and are a bit simpler. When both $k$ and $l$ are negative, the above arcs are homotopic to $(d^{-1}c^{-1})^{|k+1|}d^{-1}$ and $(f^{-1}e^{-1})^{|l+1|}e^{-1}$ respectively. The last two arcs are geodesic arcs with common end points that are homotopic as arcs with common end points. 

Assume for now that both arcs $d$ and $e$ are not the constant arcs. In this case the first geodesic arc at its starting point has the velocity vector proportional to the velocity vector of $d^{-1}$ at its start. The second geodesic arc at its starting point has the velocity vector proportional to the velocity vector of $e^{-1}$ at its start. Since $b_1$ and $b_2$ intersect transversally, the velocity vectors of the above geodesic arcs are linearly independent and hence they are different geodesic arcs. However for hyperbolic surfaces the two different geodesic segments with common ends are not homotopic, see~\cite[Theorem 1.5.3]{Buser}.

Assume that the arc $d$ is a constant arc while the arc $e$ is not a constant arc. In this case the first geodesic arc at its starting point has the velocity vector proportional to the velocity vector of $c^{-1}$ at its start. The second geodesic arc at its starting point has the velocity vector proportional to the velocity vector of $e^{-1}$ at its start. Since $b_1$ and $b_2$ intersect transversally, we use the argument similar to the one above to get the contradiction.

The case where $e$ is a constant arc while $d$ is not a constant arc is considered similarly. Note that 
we do not have to consider the case where both arcs $d$ and $e$ are constant, since in this case the two chord corresponding to $p_1$ and $p_2$ will coincide, contradicting our assumptions.

 The case described in Figure~\ref{chorddiagrams.fig}.B is the same as case \ref{chorddiagrams.fig}.A as discussed above.

\subsection{The proof of Theorem~\ref{maintheorem} in the case where $\alpha_1$ and $\alpha_2$ do contain powers of the same loop.}\label{powersofthesameelement}

In this case let $b$ be the unique (up to reparameterization) closed geodesic loop such that $\alpha_1$ contains $b^{k_1}$ and $\alpha_2$ contains $b^{k_2}.$ 

\begin{figure}
	[htbp] 
	\begin{center}

		\epsfxsize 11cm \hepsffile{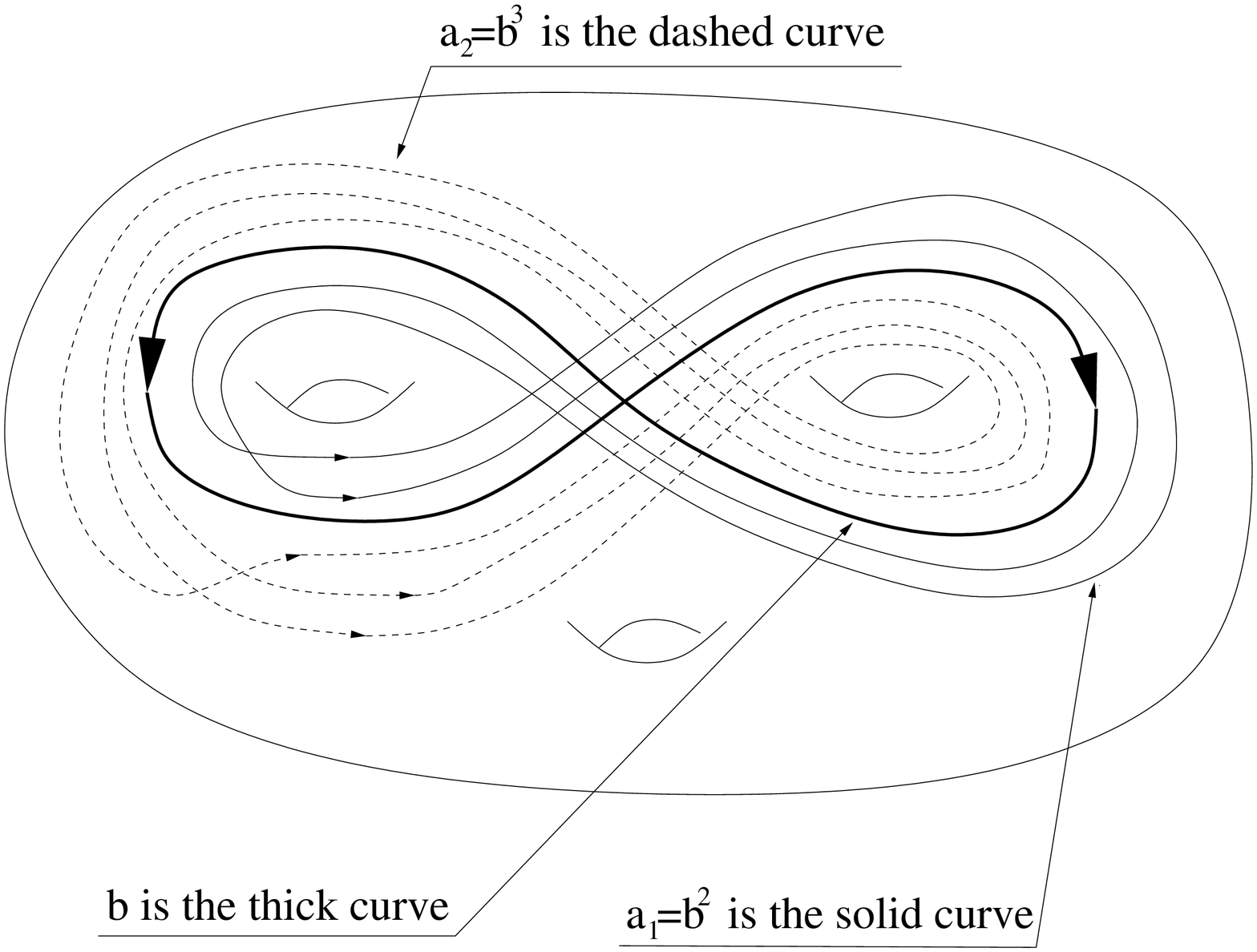} 

	\end{center}
	\caption{The curves $b$, $a_1=b^2$ and $a_2=b^3$ and the intersection points of $a_1$ and $a_2.$}\label{parallelcurves.fig} 
\end{figure}

Choose a representative $a_1$ of $\alpha_1$ that is located locally to the left of $b$, and that travels 
$k_1$ times along $b$. (Here if $k_1<0$ then this representative is oriented in the direction opposite from the direction of $b.$) Choose a representative $a_2$ of $\alpha_2$ that is located locally to the right of $b$, and that travels $k_2$ times along $b$. See Figure~\ref{parallelcurves.fig}. By changing the orientation of $b$ if necessary, {\bf we can assume that $k_1>0$.\/}

It is easy to see that if $p$ is a self intersection point of $b$ subdividing $b$ into two arcs $h_{1,p}, h_{2,p}$ then the are $2|k_1 k_2|$ intersection points of $a_1$ and $a_2$ located in a small neighborhood of $p$. These intersection points can be partitioned into two sets of $|k_1k_2|$ points each, such that the intersection points in one set have one sign and the points in the other set have the opposite sign. The total input into the bracket of these terms is: $$\pm k_1k_2\bigl((h_{1,p}h_{2,p})^{k_1}\bullet_p (h_{2,p}h_{1,p})^{k_2}\bigr)-\pm k_1k_2\bigl((h_{2,p}h_{1,p})^{k_1}\bullet_p (h_{1,p}h_{2,p})^{k_2}\bigr).$$

We have three cases to consider. First, we show a positive term and negative term corresponding to the same self-intersection point of $b$ cannot cancel. Then we show terms corresponding to different self-intersection points of $b$ cannot cancel, and this case has two subcases which depend on the Gauss diagram of $b$.

{\bf The case where two terms correspond to the same self-intersection point of $b$.} Let $c=h_{1,p}$ and $d=h_{2,p}$ be the two loops based at $p$ in the image of $b$, ordered according to the positively oriented frame formed by the tangent vectors to $b$ at $p$. Then two terms of the bracket of opposite signs are of the form $(cd)^{k_1}\bullet_p (dc)^{k_2}$ and $(dc)^{k_1}\bullet_p(cd)^{k_2}$. We suppose these terms cancel. After conjugating the first of these terms by the class of the loop $d$, we have the following equations: 
\begin{equation}
	\mu (dc)^{k_1} \mu^{-1}=(dc)^{k_1}, 
\end{equation}
and 
\begin{equation}
	\label{sameselfint} \mu d(dc)^{k_2}d^{-1}\mu^{-1}=(cd)^{k_2}. 
\end{equation}

Since $(dc)\in\pi_1(F,p)$ generates the centralizer of $(dc)^{k_1}$, then $\mu=(dc)^n$ for 
$n\in \mathbb{Z}$. Substituting this into Equation \eqref{sameselfint}, we have 
\begin{equation}
	(dc)^nd(dc)^{k_2}=(cd)^{k_2}(dc)^nd, 
\end{equation}
which implies that $(dc)^ndc^{-1}(cd)^{k_2}c=(cd)^{k_2}(dc)^nd$. Therefore $(dc)^ndc^{-1}$ and $(cd)^{k_2}$ are commuting elements of $\pi_1(F,p)$. Since $(cd)$ generates the centralizer of $(cd)^{k_2}$, we have that $(dc)^ndc^{-1}=(cd)^m$ for some $m\in \mathbb{Z}$. Hence the geodesics $(dc)^nd$ and $(cd)^mc$ are homotopic as paths with common ends. This can only happen if either $c$ or $d$ is trivial, or if one of $c$ or $d$ is a power of the other, which is impossible, because $p$ is a transverse self-intersection point of the geodesic $b$.

{\bf The cases where the two terms correspond to different self-intersection points of $b$.}

The terms corresponding to two distinct self intersection points $p_1$ and $p_2$ of $b$ are 
$$\epsilon k_1k_2\bigl((h_{1,p_1}h_{2,p_1})^{k_1}\bullet_{p_1} (h_{2,p_1}h_{1,p_1})^{k_2}\bigr)-\epsilon k_1k_2\bigl((h_{2,p_1}h_{1,p_1})^{k_1}\bullet_{p_1} (h_{1,p_1}h_{2,p_1})^{k_2}\bigr),$$
and
$$\delta k_1k_2\bigl((h_{1,p_2}h_{2,p_2})^{k_1}\bullet_{p_2} (h_{2,p_2}h_{1,p_2})^{k_2}\bigr)-\delta k_1k_2\bigl((h_{2,p_2}h_{1,p_2})^{k_1}\bullet_{p_2} (h_{1,p_2}h_{2,p_2})^{k_2}\bigr),$$
for $\epsilon,\delta=\pm 1$, where, for each $j=1,2$, $h_{1,p_j}$ and $h_{2,p_j}$ are the two loops in $b$ based at $p_j$.  For this part of the proof, it does not matter which loop we label $h_{1,p_j}$ and which loop we label $h_{2,p_j}$.

By symmetry, it suffices to show that the geometric chord diagrams $(h_{1,{p_1}}h_{2,{p_1}})^{k_1}\bullet _{p_1} (h_{2,{p_1}}h_{1,{p_1}})^{k_2}$ and $(h_{2,{p_2}}h_{1,{p_2}})^{k_1}\bullet _{p_2}(h_{1,{p_2}}h_{2,{p_2}})^{k_2}$ are not homotopic, {\it and} that the geometric chord diagrams $(h_{1,{p_1}}h_{2,{p_1}})^{k_1}\bullet _{p_1} (h_{2,{p_1}}h_{1,{p_1}})^{k_2}$ and $(h_{1,{p_2}}h_{2,{p_2}})^{k_1}\bullet _{p_2}(h_{2,{p_2}}h_{1,{p_2}})^{k_2}$ are not homotopic. We argue by contradiction and assume that one of these pairs is formed by homotopic chord diagrams.

\begin{figure}
	[htbp] 
	\begin{center}

		\epsfxsize 10cm \hepsffile{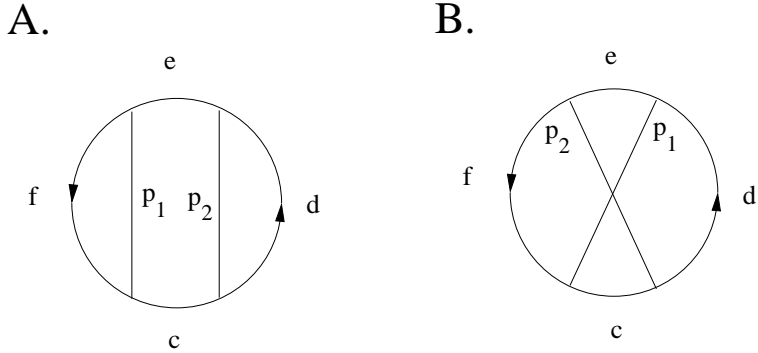} 

	\end{center}
	\caption{Two possible mutual positions for chords corresponding to self-intersection points of $b$}\label{chorddiagramsself.fig} 
\end{figure}

There are two possible cases for the mutual position of chords of $b$ corresponding to $p_1$ and $p_2.$ They are shown in Figure~\ref{chorddiagramsself.fig}.A and Figure~\ref{chorddiagramsself.fig}.B respectively. We consider these two cases separately.

As we already stated, for our proof it does not matter which of the two loops 
corresponding to the point $p_j$ 
we label $h_{1,p_j}$ and 
which loop we label $h_{2,p_j}$. So for the case described on Figure~\ref{chorddiagramsself.fig}.A
we assume that $h_{1,p_1}=cde,$ $h_{2, p_1}=f$, $h_{1, p_2}=efc$, and $h_{2, p_2}=d.$ For the case
described on Figure~\ref{chorddiagramsself.fig}.B
we assume that $h_{1,p_1}=cd,$ $h_{2, p_1}=ef$, $h_{1, p_2}=fc$, and $h_{2, p_2}=de.$

{\bf The case described in Figure~\ref{chorddiagramsself.fig}.A: }

{\bf First, we show the pair of chord diagrams $(h_{1,{p_1}}h_{2,{p_1}})^{k_1}\bullet _{p_1} (h_{2,{p_1}}h_{1,{p_1}})^{k_2}$ and $(h_{2,{p_2}}h_{1,{p_2}})^{k_1}\bullet _{p_2}(h_{1,{p_2}}h_{2,{p_2}})^{k_2}$ cannot cancel.} Suppose these terms do cancel.  In this case we get that 
\begin{equation}
	(h_{1,{p_1}}h_{2,{p_1}})^{k_1}\bullet _{p_1} (h_{2,{p_1}}h_{1,{p_1}})^{k_2}= (cdef)^{k_1}\bullet _{p_1}(fcde)^{k_2} 
\end{equation}
and 
\begin{equation}
	(h_{2,{p_2}}h_{1,{p_2}})^{k_1}\bullet _{p_2}(h_{1,{p_2}}h_{2,{p_2}})^{k_2}= (defc)^{k_1}\bullet _{p_2} (efcd)^{k_2} 
\end{equation}
We deform the second chord diagram, so that its chord is mapped to $p_1$:
\begin{equation}
	(defc)^{k_1}\bullet _{p_2} (efcd)^{k_2} = c(defc)^{k_1}c^{-1}\bullet _{p_1} c(efcd)^{k_2}c^{-1} = (cdef)^{k_1}\bullet _{p_1}c(efcd)^{k_2}c^{-1}. 
\end{equation}
Now both chord diagrams can be considered as pairs of loops based at $p_1$; see Equation~\eqref{descr}.  Since we assumed that 
\begin{equation}
	(h_{1,{p_1}}h_{2,{p_1}})^{k_1}\bullet _{p_1} (h_{2,{p_1}}h_{1,{p_1}})^{k_2} = (h_{2,{p_2}}h_{1,{p_2}})^{k_1}\bullet _{p_2}(h_{1,{p_2}}h_{2,{p_2}})^{k_2} 
\end{equation}

and $\alpha_1\neq \alpha_2$, we get that there exists $\mu\in \pi_1(F, p_1)$ such that 
\begin{equation}
	\label{caseA1} \mu(cdef)^{k_1}\mu^{-1}=(cdef)^{k_1}\in \pi_1(F, p_1) 
\end{equation}
and 
\begin{equation}
	\label{caseA2} \mu c (efcd)^{k_2} c^{-1}\mu^{-1}=(fcde)^{k_2}\in \pi_1(F, p_1). 
\end{equation}

Since the centralizer of $(cdef)^{k_1}\in \pi_1(F, p_1)$ is the free abelian group generated by $(cdef)\in\pi_1(F, p_1)$, we have that $\mu=(cdef)^k\in \pi_1(F, p_1)$ for some $k\in\Z.$ Substitute this into equation~\eqref{caseA2} to get 
\begin{equation}
	\label{splitk2} (cdef)^k c (efcd)^{k_2} c^{-1}(cdef)^{-k}=(fcde)^{k_2}\in \pi_1(F, p_1). 
\end{equation}

It follows from the previous equation that 
\begin{equation}
	(cdef)^k ce (fcde)^{k_2} e^{-1}c^{-1}(cdef)^{-k}=(fcde)^{k_2}\in \pi_1(F, p_1). 
\end{equation}
Note that this holds even if $k_2$ is negative.

Thus $(cdef)^kce(fcde)^{k_2}=(fcde)^{k_2}(cdef)^kce\in \pi_1(F, p_1)$ and the loops $(cdef)^kce$ and $(fcde)^{k_2}$ commute. Since the centralizer of $(fcde)^{k_2}$ is generated by $fcde$ we have that $(cdef)^k ce=(fcde)^l\in\pi_1(F, p_1)$ for some $l\in\Z$. 

We consider the case when both $k$ and $l$ are negative, the other cases are considered in a similar fashion and are a bit simpler. We have that $(f^{-1}e^{-1}d^{-1}c^{-1})^{|k|}c= (f^{-1}e^{-1}d^{-1}c^{-1})^{|k|-1}f^{-1}e^{-1}d^{-1}$ and $(e^{-1}d^{-1}c^{-1}f^{-1})^{|l|}e^{-1}$ are geodesic arcs with the same start and end points and they are homotopic as paths with fixed starting and ending points.

Note that $f$ and $d$ are never constant loops because they correspond to transverse self-intersection 
points of the geodesic $b$. On the other hand, since we assumed that the self intersection points 
$p_1, p_2$ are different, exactly one of the arcs $e$ or $c$ could be trivial.  Assume for now that both arcs $c$ and $e$ are not the constant arcs. In this case the first geodesic arc at its starting point has the velocity vector proportional to the velocity vector of $f^{-1}$ at its start. The second geodesic arc at its starting point has the velocity vector proportional to the velocity vector of $e^{-1}$ at its start. Since the branches of $b$ intersect transversally at the self intersection points the velocity vectors of the above geodesic arcs are linearly independent and hence they are different geodesic arcs. However for hyperbolic surfaces the two different geodesic segements with common ends are not homotopic, see~\cite[Theorem 1.5.3]{Buser}.  If $e$ is a constant arc, while $d$ is nontrivial, then one can check that the velocity vectors of the above geodesic arcs at the end points are still linearly independent and we again get a contradiction. The case where $d$ is a trivial arc, while $e$ is nontrivial is considered similarly.

{\bf Next, we show the pair of chord diagrams $(h_{1,{p_1}}h_{2,{p_1}})^{k_1}\bullet _{p_1} (h_{2,{p_1}}h_{1,{p_1}})^{k_2}$ and $(h_{1,{p_2}}h_{2,{p_2}})^{k_1}\bullet _{p_2}(h_{2,{p_2}}h_{1,{p_2}})^{k_2}$ cannot cancel.} Suppose these terms do cancel.  In this case we get that 
\begin{equation}
	(h_{1,{p_1}}h_{2,{p_1}})^{k_1}\bullet _{p_1} (h_{2,{p_1}}h_{1,{p_1}})^{k_2}= (cdef)^{k_1}\bullet _{p_1}(fcde)^{k_2} 
\end{equation}
and 
\begin{equation}
	(h_{1,{p_2}}h_{2,{p_2}})^{k_1}\bullet _{p_2}(h_{2,{p_2}}h_{1,{p_2}})^{k_2}= (efcd)^{k_1}\bullet _{p_2} (defc)^{k_2} 
\end{equation}
Now we deform the second chord diagram so that its chord is mapped to $p_1$:
\begin{equation}
\begin{split}
	  (efcd)^{k_1}\bullet _{p_2} (defc)^{k_2}= cd(efcd)^{k_1}d^{-1}c^{-1}\bullet _{p_1} cd(defc)^{k_2}d^{-1}c^{-1}\\=(cdef)^{k_1}\bullet _{p_1} cd(defc)^{k_2}d^{-1}c^{-1}.
\end{split}
\end{equation}
Now both chord diagrams can be considered as pairs of loops based at $p_1$; see Equation~\eqref{descr}.  Since we assumed that 
\begin{equation}
\begin{split}
	(h_{1,{p_1}}h_{2,{p_1}})^{k_1}\bullet _{p_1} (h_{2,{p_1}}h_{1,{p_1}})^{k_2} = (h_{1,{p_2}}h_{2,{p_2}})^{k_1}\bullet _{p_2}(h_{2,{p_2}}h_{1,{p_2}})^{k_2} 
\end{split}
\end{equation}

and $\alpha_1\neq \alpha_2$, we get that there exists $\mu\in \pi_1(F, p_1)$ such that 
\begin{equation}
	\label{diagram4Acase2eqn1} \mu(cdef)^{k_1}\mu^{-1}=(cdef)^{k_1}\in \pi_1(F, p_1) 
\end{equation}
and 
\begin{equation}
	\label{diagram4Acase2eqn2} \mu cd(defc)^{k_2}d^{-1}c^{-1}\mu^{-1}=(fcde)^{k_2}\in \pi_1(F, p_1). 
\end{equation}
Since the centralizer of $(cdef)^{k_1}$ is the free abelian group generated by $cdef\in \pi_1(F,p_1)$, we get that $\mu=(cdef)^k$ for some $k\in \Z$.  We substitute this into Equation \eqref{diagram4Acase2eqn2} to get 
$$ (cdef)^k cd(defc)^{k_2}d^{-1}c^{-1}(cdef)^{-k}=(fcde)^{k_2}\in \pi_1(F, p_1).$$
We can rewrite this equation to get 
\begin{equation}
c^{-1}f^{-1}(cdef)^k cd(defc)^{k_2}=(defc)^{k_2}c^{-1}f^{-1}(cdef)^kcd\in \pi_1(F, p_2).
\end{equation}
  Therefore $c^{-1}f^{-1}(cdef)^kcd$ commutes with $(defc)^{k_2}$ in $\pi_1(F, p_2)$, so 
\begin{equation}c^{-1}f^{-1}(cdef)^kcd=(defc)^l\in\pi_1(F, p_2)
\end{equation}
 for some $l\in\Z$.  Hence the paths $(cdef)^kcd$ and $fc(defc)^l$ are homotopic as paths with fixed ends. The arcs $f$ and $d$ cannot be trivial, but one of the arcs $d$ or $e$ might be trivial. 

Now we consider the cases of various possible signs of $k$ and $l$. 

For example if both $k$ and $l$ negative, we get that the paths 
$(f^{-1}e^{-1}d^{-1}c^{1})^{|k|}cd=
(f^{-1}e^{-1}d^{-1}c^{1})^{|k|-1}f^{-1}e^{-1}$ and $fc(c^{-1}f^{-1}e^{-1}d^{-1})^{|l|}=(e^{-1}d^{-1}c^{-1}
f^{-1})^{|l|-1}e^{-1}d^{-1}$ are homotopic as geodesic arcs with fixed end points. We observe that the 
velocity vectors of these arcs at the starting points are linearly independent, hence these are different 
geodesic arcs and they can not be homotopic by~\cite[Theorem 1.5.3]{Buser}.

The cases of other signs for $k$ and $l$ are considered similarly.

{\bf The case described in Figure~\ref{chorddiagramsself.fig}.B.}

{\bf First we show the pair chord  diagrams $(h_{1,{p_1}}h_{2,{p_1}})^{k_1}\bullet _{p_1} (h_{2,{p_1}}h_{1,{p_1}})^{k_2}$ and $(h_{2,{p_2}}h_{1,{p_2}})^{k_1}\bullet _{p_2}(h_{1,{p_2}}h_{2,{p_2}})^{k_2}$ cannot cancel.} Suppose these terms do cancel.  In this case we get that 
\begin{equation}
	(h_{1,{p_1}}h_{2,{p_1}})^{k_1}\bullet _{p_1} (h_{2,{p_1}}h_{1,{p_1}})^{k_2} = (cdef)^{k_1}\bullet _{p_1}(efcd)^{k_2} 
\end{equation}
and 
\begin{equation}
	(h_{2,{p_2}}h_{1,{p_2}})^{k_1}\bullet _{p_2}(h_{1,{p_2}}h_{2,{p_2}})^{k_2} = (defc)^{k_1}\bullet _{p_2} (fcde)^{k_2} 
\end{equation}
Now we deform the second diagram so that its chord is mapped to $p_1$:
\begin{equation}
	(defc)^{k_1}\bullet _{p_2} (fcde)^{k_2} = c(defc)^{k_1}c^{-1}\bullet _{p_1} c(fcde)^{k_2}c^{-1} = (cdef)^{k_1}\bullet _{p_1}c(fcde)^{k_2}c^{-1} 
\end{equation}
Now both chord diagrams can be considered as pairs of loops based at $p_1$; see Equation~\eqref{descr}.  Since we assumed that 
\begin{equation}
	(h_{1,{p_1}}h_{2,{p_1}})^{k_1}\bullet _{p_1} (h_{2,{p_1}}h_{1,{p_1}})^{k_2} = (h_{2,{p_2}}h_{1,{p_2}})^{k_1}\bullet _{p_2}(h_{1,{p_2}}h_{2,{p_2}})^{k_2} 
\end{equation}

and $\alpha_1\neq \alpha_2$, we get that there exists $\mu\in \pi_1(F, p_1)$ such that 
\begin{equation}
	\label{caseB1} \mu(cdef)^{k_1}\mu^{-1}=(cdef)^{k_1}\in \pi_1(F, p_1) 
\end{equation}
and 
\begin{equation}
	\label{caseB2} \mu c (fcde)^{k_2} c^{-1}\mu^{-1}=(efcd)^{k_2}\in \pi_1(F, p_1). 
\end{equation}

Since the centralizer of $(cdef)^{k_1}\in \pi_1(F, p_1)$ is the free abelian group generated by $(cdef)\in\pi_1(F, p_1)$, we have that $\mu=(cdef)^k\in \pi_1(F, p_1)$ for some $k\in\Z.$ Substitute this into equation~\eqref{caseA2} to get 
\begin{equation}
	\label{secondsplitk2} (cdef)^k c (fcde)^{k_2} c^{-1}(cdef)^{-k}=(efcd)^{k_2}\in \pi_1(F, p_1). 
\end{equation}

From the previous equation we have 
\begin{equation}
	(cdef)^k ce^{-1} (efcd)^{k_2} ec^{-1}(cdef)^{-k}=(efcd)^{k_2}\in \pi_1(F, p_1). 
\end{equation}

Thus $(cdef)^kce^{-1}(efcd)^{k_2}=(efcd)^{k_2}(cdef)^kce^{-1}\in \pi_1(F, p_2)$ and the loops $(cdef)^kce^{-1}$ and $(efcd)^{k_2}$ commute. Since the centralizer of $(efcd)^{k_2}$ is generated by $efcd$ we have that $(cdef)^k ce^{-1}=(efcd)^l\in\pi_1(F, p_2)$ for some $l\in\Z$. 

Therefore the arcs $(cdef)^kc$ and $(efcd)^le$ are homotopic as arcs with fixed end points. Similarly to the proof in the case described on Figure~\ref{chorddiagramsself.fig}.A, we consider possible signs for $k$ and $l$. For all the possible pairs of signs, the expressions for the arcs above simplify to expressions for homotopic geodesic arcs. Considering the velocity vectors at the starting points, we see that these homotopic geodesic arcs are distinct.  This contradicts to~\cite[Theorem 1.5.3]{Buser}.

{\bf Finally, we suppose the chord diagrams $(h_{1,{p_1}}h_{2,{p_1}})^{k_1}\bullet _{p_1} (h_{2,{p_1}}h_{1,{p_1}})^{k_2}$ and $(h_{1,{p_2}}h_{2,{p_2}})^{k_1}\bullet _{p_2}(h_{2,{p_2}}h_{1,{p_2}})^{k_2}$ cancel.} 
In this case we get that 
\begin{equation}
	(h_{1,{p_1}}h_{2,{p_1}})^{k_1}\bullet _{p_1} (h_{2,{p_1}}h_{1,{p_1}})^{k_2}= (cdef)^{k_1}\bullet _{p_1}(efcd)^{k_2} 
\end{equation}
and 
\begin{equation}
	(h_{1,{p_2}}h_{2,{p_2}})^{k_1}\bullet _{p_2}(h_{2,{p_2}}h_{1,{p_2}})^{k_2}= (fcde)^{k_1}\bullet _{p_2} (defc)^{k_2} 
\end{equation}
Now we deform the first chord diagram so that its chord is mapped to $p_2$:
\begin{equation}
	  (cdef)^{k_1}\bullet _{p_1}(efcd)^{k_2}=f(cdef)^{k_1}f^{-1}\bullet _{p_2}f(efcd)^{k_2}f^{-1}=(fcde)^{k_1}\bullet_{p_2}f(efcd)^{k_2}f^{-1}
\end{equation}
Now both chord diagrams can be considered as pairs of loops based at $p_2$; see Equation~\eqref{descr}.
Since we assumed that 
\begin{equation}
	(h_{1,{p_1}}h_{2,{p_1}})^{k_1}\bullet _{p_1} (h_{2,{p_1}}h_{1,{p_1}})^{k_2} = (h_{1,{p_2}}h_{2,{p_2}})^{k_1}\bullet _{p_2}(h_{2,{p_2}}h_{1,{p_2}})^{k_2} 
\end{equation}

and $\alpha_1\neq \alpha_2$, we get that there exists $\mu\in \pi_1(F, p_2)$ such that 
\begin{equation}
	\label{diagram4Bcase2eqn1} \mu(fcde)^{k_1}\mu^{-1}=(fcde)^{k_1}\in \pi_1(F, p_2) 
\end{equation}
and 
\begin{equation}
	\label{diagram4Bcase2eqn2} \mu f(efcd)^{k_2}f^{-1}\mu^{-1}=(defc)^{k_2}\in \pi_1(F, p_2). 
\end{equation}
Since the centralizer of $(fcde)^{k_1}$ is the free abelian group generated by $fcde\in \pi_1(F,p_2)$, we get that $\mu=(fcde)^k$ for some $k\in \Z$.  We substitute this into Equation \eqref{diagram4Bcase2eqn2} to get 
$$(fcde)^k f(efcd)^{k_2}f^{-1}(fcde)^{-k}=(defc)^{k_2}\in \pi_1(F, p_2).$$
We can rewrite this equation to get 
\begin{equation}
d^{-1}(fcde)^k f(efcd)^{k_2}=(efcd)^{k_2}d^{-1}(fcde)^kf\in \pi_1(F, p_2).  
\end{equation}
Therefore $d^{-1}(fcde)^k f$ commutes with $(efcd)^{k_2}$ in $\pi_1(F, p_2)$, so $d^{-1}(fcde)^k f=(efcd)^l\in \pi_1(F, p_2)$ for some $l\in\Z$.  Thus $(fcde)^kf$ and $d(efcd)^l$ are homotopic as paths with fixed end points. 
Now the end of the proof is similar to how we finished the proof in the case on Figure~\ref{chorddiagramsself.fig}.B; equality of chord  diagrams $(h_{1,{p_1}}h_{2,{p_1}})^{k_1}\bullet _{p_1} (h_{2,{p_1}}h_{1,{p_1}})^{k_2}$ and $(h_{2,{p_2}}h_{1,{p_2}})^{k_1}\bullet _{p_2}(h_{1,{p_2}}h_{2,{p_2}})^{k_2}$.

This finishes the proof of Theorem~\ref{maintheorem} for all compact surfaces $F$ with or without boundary. 

{\bf Now we consider the case where $F$ is any noncompact surface.\/} Note that $F$ may be much more complicated than a compact surface with the boundary deleted.

Given two different free homotopy classes $\alpha_1$ and $\alpha_2$ of curves on $F$, we choose generic curves $b_1\in \alpha_1$ and $b_2\in \alpha_2$ on $F$ such that $i(b_1, b_2)=\#(\alpha_1, \alpha_2).$ We argue by contradiction and assume that $\term(\{\alpha_1, \alpha_2\})$ is strictly less than $\#(\alpha_1, \alpha_2).$ Then there are two different intersection points $p, p'$ of $b_1$ and $b_2$ such that the signs of the intersection points are opposite but the corresponding two chord diagrams, each one been of the form two circles connected by a chord, are homotopic.

Let $\widetilde F\subset F$ be a compact surface with boundary that contains $b_1, b_2$ and the image of the homotopy between the two chord diagrams described above. Put $\tilde \alpha_i$ to be the free homotopy class of $b_i, i=1,2$ on the surface $\tilde F.$

Clearly $\#(\tilde \alpha_1, \tilde \alpha_2)$ computed with respect to the surface $\widetilde F$ equals to $\#(\alpha_1, \alpha_2)$ computed with respect to the surface $F.$ Moreover $i(b_1, b_2)=\#(\tilde \alpha_1, \tilde \alpha_2)$. By our choice of $\widetilde F$, the two chord diagrams coming from the intersection points $p, p'$ of $b_1, b_2$ are homotopic within $\widetilde F$ and the signs of these intersection points are different. Thus for the compact surface $\widetilde F$ we have $\term(\{\tilde \alpha_1, \tilde \alpha_2\})<i(b_1, b_2)=\#(\tilde \alpha_1, \tilde \alpha_2).$ This contradicts the statement of our theorem that is already proved for all compact surfaces $\widetilde F.$ 

This finishes the proof of Theorem~\ref{maintheorem} for all surfaces. \qed

\section{Proof of Theorem \ref{mainthmcor}}
Let $\alpha$ be a nontrivial free homotopy class of loops on the oriented surface $F$ and let $n$ be the largest integer such that $\alpha=\beta^n$ for some $\beta\in \pi_1(F)$.  Let $p$ and $q$ be distinct nonzero integers.  Note that we allow the case where $p$ and $q$ have the same absolute value.  We may assume, as in the previous section, that $F$ is compact, with or without boundary.  If $F=A^2, T^2$ or $S^2$ then $\alpha$ must be a power of a simple class, and hence $\#(\alpha)=n-1$; see Hass-Scott \cite{HassScottShortening} or \cite{TuraevViro}. Now suppose $F\neq A^2, T^2,$ or $S^2$. Let $b$ be the unique geodesic representative of $\beta$.   Let $r_k$ be a degree $k$ curve on the annulus with $k-1$ self-intersection points (which, by \cite{HassScottShortening}, is the unique such curve up to isotopy), and let $b_k=f\circ r_k$ be a curve on $F$ given by composing $r_k$ with a map $f:A^2\rightarrow F$ that maps the annulus into a thin regular neighborhood of $b$.  Put $a=b_n$.  The number of self-intersection points $i(b_k)$ of $b_k$ is $i(b)k^2+k-1$, and in particular $i(a)=i(b)n^2+n-1$, where $i(b)$ is the number of self-intersection points of $b$.   Hass and Scott \cite{HassScottShortening} show that $i(b_k)=\#(\beta^k)$, and in particular, that $i(a)=\#(\alpha)$.
Now let $a_1$ be the curve obtained by shifting $b_{np}$ slightly so that it is located locally to the left of $b$.  Let $a_2$ be the curve obtained by shifting $b_{nq}$ slightly so that it is located locally to the right of $b$.  The curves $a_1$ and $a_2$ are simply the curves described in Subsection \ref{powersofthesameelement} for the case where $k_1=np$ and $k_2=nq$.  By the argument in Subsection \ref{powersofthesameelement}, we have $i(a_1,a_2)=2n^2|pq|i(b)$, and $\term(\{\alpha^p,\alpha^{q}\})=2n^2|pq|i(b)$.   Hence $\#(\alpha)=\frac{1}{2|pq|}\term(\{\alpha^p,\alpha^{q}\})+n-1$.\qed

\section{Computing the Minimal Intersection Number Explicitly Using $\{ \cdot, \cdot \}$}\label{example} In this section we show how to compute the minimal number of intersection points of loops in two given free loop homotopy classes $\alpha_1, \alpha_2$ using the Andersen-Mattes-Reshetkhin Poisson bracket $\{\alpha_1, \alpha_2 \}$. The pair $\alpha_1, \alpha_2$ of homotopy classes we use was found by Chas~\cite[Example 5.5]{ChasCombinatorial}. Chas showed that in this example the Goldman Lie bracket vanishes, while $\#(\alpha_1, \alpha_2)\neq 0.$ In particular, in this example the Goldman Lie bracket does not allow one to compute the minimal number of intersection points of loops in two given free homotopy classes. She used an algorithm~\cite[Theorem 3.13, Remark 3.14]{ChasCombinatorial} for finding the minimal number of intersection points on a surface with boundary to show that the minimal number of intersection points in this example is $2.$ ( Chas used an algorithm similar to the one constructed by M.~Cohen and M.~Lustig~\cite{CohenLustig}.) We get the same answer using the Andersen-Mattes-Reshetikhin Poisson bracket $\{\alpha_1, \alpha_2 \}.$ 

The problem of finding the minimal number of intersection points is trivial if $F=S^2, A \text{ (annulus)}$ or if one of $\alpha_1, \alpha_2$ is the class $1$ of the constant loop.

For $F=T^2$ the number $\#(\alpha_1, \alpha_2)$ is easy to compute, since it equals the absolute value of the intersection number of the classes in $H_1(T^2)$ realized by $\alpha_1, \alpha_2;$ see, for example, the proof of Theorem~\ref{maintheorem} for $F=T^2.$ Of course for $F=T^2$ we also have $\#(\alpha_1, \alpha_2)=\term(\{\alpha_1, \alpha_2\})$, again by Theorem~\ref{maintheorem}. 

Let us outline the general procedure of how one uses $\{\alpha_1, \alpha_2 \}$ to find $\#(\alpha_1, \alpha_2)$, for $F\neq A, S^2, T^2$ and distinct $\alpha_1, \alpha_2\neq 1.$ 

By Theorem~\ref{maintheorem} $\term(\{\alpha_1, \alpha_2\})=\#(\alpha_1, \alpha_2).$ 

For $h_1, h_2\in \pi_1(F, p)$ we denote by $h_1\bullet_p h_2$ the homotopy class of a chord diagram as described 
in Subsection~\ref{bulletanditsproperties}.

  Assume that the terms $+1(h_1\bullet_p h_2)$ and $-1(\tilde h_1\bullet_p \tilde h_2)$ cancel, i.e. that 
$h_1\bullet_p h_2=\tilde h_1\bullet_p \tilde h_2.$ We have 
assumed that the free homotopy classes $\alpha_1, \alpha_2$ are different. Thus by Equation~\eqref{descr}
we may assume that $\langle h_1\rangle=\langle \tilde h_1\rangle =\alpha_1$,
while $\langle h_2 \rangle=\langle \tilde h_2\rangle=\alpha_2$.

Since $\langle h_1\rangle=\langle \tilde h_1\rangle =\alpha_1$
 we choose $g\in \pi_1(F, p)$ such that $g \tilde h_1 g^{-1}=h_1$. We get that
\begin{equation}
	\tilde h_1\bullet_p \tilde h_2=g\tilde h_1g^{-1}\bullet_p g\tilde h_2g^{-1}=h_1\bullet_p g\tilde h_2g^{-1}=h_1\bullet_p h_2. 
\end{equation}

Thus by equation~\eqref{descr} there exists $\tilde s\in \pi_1(F, p)$ such that 
\begin{equation}
	\label{usefulequation} \tilde sh_1 \tilde s^{-1}=h_1 \text{ and }\tilde s(g \tilde h_2 g^{-1})\tilde s^{-1}=h_2. 
\end{equation}

The nontrivial abelian subgroups of $\pi_1(F)$ are infinite cyclic, and each $\beta\neq 1\in \pi_1(F,p)$ is contained in a unique, maximal infinite cyclic group. For closed $F\neq S^2, T^2$ this follows from the Preissman Theorem~\cite[Theorem 3.2 page 260]{docarmo}, since such surfaces admit a hyperbolic Riemannian metric. For general noncompact surfaces (and hence for compact surfaces with boundary) this is true, since $\pi_1(F)$ is free ~\cite[pages 142-144]{Stillwell},~\cite{Johannsson}.

Let $s$ be a generator of the maximal infinite cyclic group containing $h_1.$ Equation~\eqref{usefulequation} says that $\tilde s$ and $h_1$ commute, so $\tilde s=s^i\in \pi_1(F,p)$ for some $i\in \Z.$ Now one checks if there exists some $i\in \Z$ such that $s^i(g \tilde h_2 g^{-1}) s^{-i}=h_2.$ If such $i$ exists, then the terms $+1(h_1\bullet_p h_2)$ and $-1(\tilde h_1\bullet_p \tilde h_2)$ cancel; if such $i$ does not exist, then the terms do not cancel. 

\begin{figure}
	[htbp] 
	\begin{center}

		\epsfxsize 10cm \hepsffile{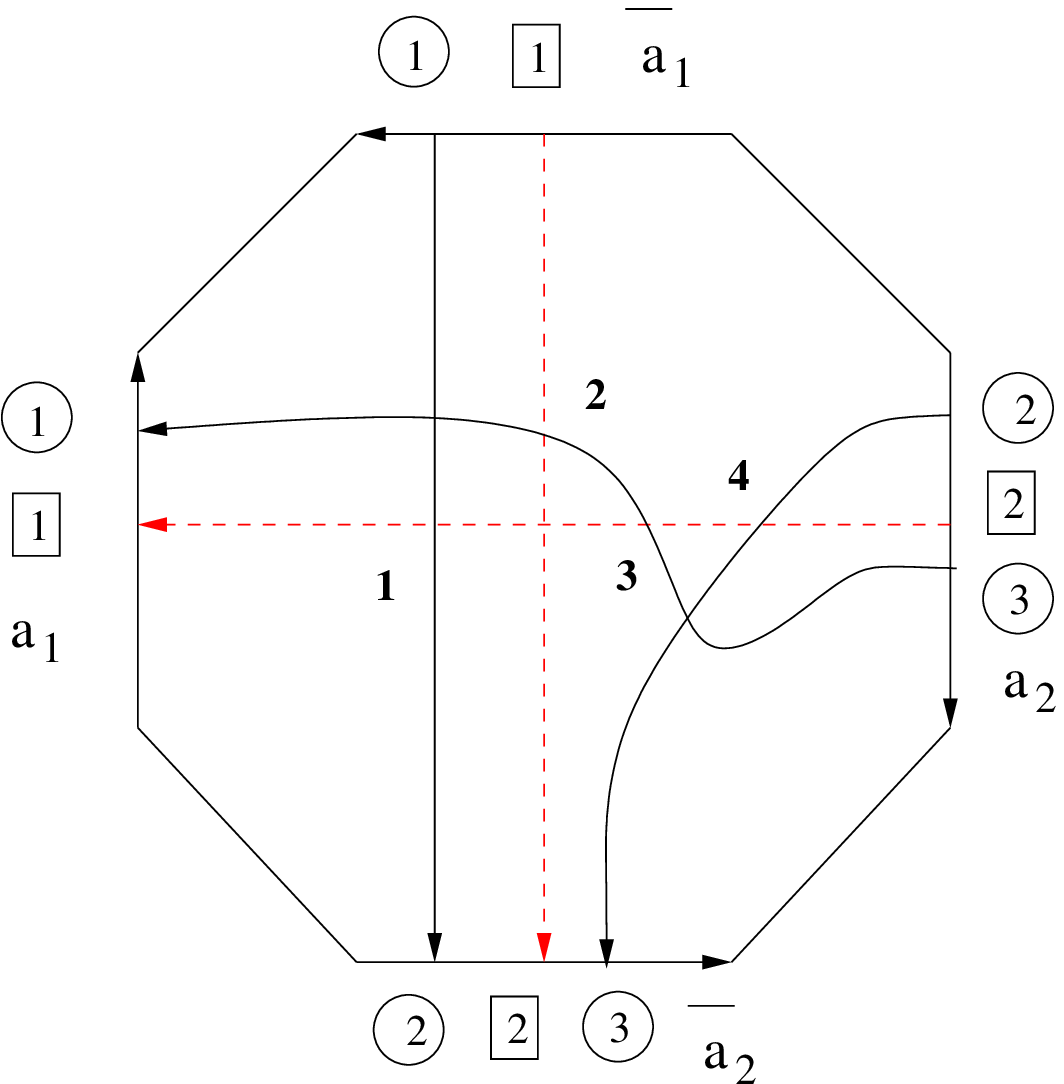} 

	\end{center}
	\caption{}\label{interexample0.fig} 
\end{figure}

As an example, consider the oriented surface $F$ with boundary obtained from the $8$-gon shown in Figure~\ref{interexample0.fig} by gluing the side $a_1$ to $\ov a_1$ and the side $a_2$ to $\ov a_2$, in such a way that the points corresponding to the arrow heads on $a_i$ and on $\ov a_i$ are identified. We put $p$ to be the center of the $8$-gon. Clearly $\pi_1(F,p)$ is a free group on two generators. We denote the generators by $a_1$ and $a_2$ according to the following convention introduced in~\cite{ChasCombinatorial}: $a_i\in \pi_1(F)$ is the class of the loop that starts from $p$, goes to the side $a_i,$ passes through the image of this side on $F,$ reappears in the $8$-gon from the side denoted by $\ov a_i$ and returns back to $p.$ 

Let $b_1$ and $b_2$ be respectively the solid and the dashed curves in Figure~\ref{interexample0.fig} and $\alpha_1,\alpha_2$ be the free homotopy classes realized by them. The numbers in the circles (respectively in the boxes) in this Figure indicate the points of $b_1$ (respectively of $b_2$) on the sides of the $8$-gon that have to be identified. The numbers $1,2,3,4$ in the middle of the $8$-gon enumerate the $4$ intersection points between the loops $b_1$ and $b_2.$

From Figure~\ref{interexample0.fig} we get that 
\begin{equation}
	\label{computebracket} 
	\begin{split}
		\{\alpha_1, \alpha_2\}= a_2^{-2}a_1\bullet_p a_1a_2^{-1} - a_1a_2^{-2}\bullet_p a_2^{-1}a_1-a_1a_2^{-2}\bullet_p a_1a_2^{-1}+ a_2^{-1}a_1a_2^{-1}\bullet_p a_1a_2^{-1}, 
	\end{split}
\end{equation}
where the $4$ terms of the expression are given in the order of the enumerated intersection points of $b_1$ and $b_2.$

Conjugating both loops by $a_1a_2^{-1}$ we get that the $4$-th term of~\eqref{computebracket} is 
\begin{equation}
	(a_1a_2^{-1})(a_2^{-1}a_1a_2^{-1})(a_1a_2^{-1})^{-1}\bullet_p (a_1a_2^{-1})(a_1a_2^{-1})(a_1a_2^{-1})^{-1}= a_1a_2^{-2}\bullet_p a_1a_2^{-1}. 
\end{equation}
Thus the third and the fourth term of~\eqref{computebracket} cancel and 
\begin{equation}
	\label{simplebracket} \{\alpha_1, \alpha_2\}=a_2^{-2}a_1\bullet_p a_1a_2^{-1} - a_1a_2^{-2}\bullet_p a_2^{-1}a_1. 
\end{equation}

Let us show that the chord diagrams $a_2^{-2}a_1\bullet_p a_1a_2^{-1}$ and $a_1a_2^{-2}\bullet_p a_2^{-1}a_1$ are not homotopic. 

Conjugate by $a_2^2$ to get that $a_2^{-2}a_1\bullet_p a_1a_2^{-1}=a_2^2(a_2^{-2}a_1)a_2^{-2}\bullet_p a_2^2(a_1a_2^{-1})a_2^{-2}= a_1a_2^{-2}\bullet_p a_2^2a_1a_2^{-3}.$ If $a_2^{-2}a_1\bullet_p a_1a_2^{-1}=a_1a_2^{-2}\bullet_p a_2^2a_1a_2^{-3}$ and $a_1a_2^{-2}\bullet_p a_2^{-1}a_1$ are equal, then there exists $\tilde s\in \pi_1(F,p)$ such that 
\begin{equation}
	\label{useful2equation} \tilde s(a_1a_2^{-2})\tilde s^{-1}=a_1a_2^{-2} \text{ and }\tilde s(a_2^2a_1a_2^{-3})\tilde s^{-1} =a_2^{-1}a_1. 
\end{equation}
The class in $H_1(F)$ realized by $a_1a_2^{-2}$ is not realizable as a nontrivial integer multiple of any element other than itself. Hence the maximal infinite cyclic group containing $a_1a_2^{-2}$ is generated by $s=a_1a_2^{-2}.$

Equation~\eqref{useful2equation} says that $\tilde s$ and $a_1a_2^{-2}$ commute. As we discussed above, this implies that $\tilde s$ and $a_1a_2^{-2}$ are in the same infinite cyclic group. Hence $\tilde s=s^i=(a_1a_2^{-2})^i,$ for some $i\in \Z.$

Thus by Equation~\eqref{useful2equation} we should have $(a_1a_2^{-2})^{i}(a_2^2a_1a_2^{-3})(a_1a_2^{-2})^{-i}=a_2^{-1}a_1,$ for some $i\in \Z.$ Since $a_1$ and $a_2$ are generators of the free group $\pi_1(F),$ it is easy to see that this identity could not hold for any $i\in \Z.$ We got a contradiction and hence the chord diagrams $a_2^{-2}a_1\bullet_p a_1a_2^{-1}$ and $a_1a_2^{-2}\bullet_p a_2^{-1}a_1$ are not homotopic. 

Hence $\term(\{\alpha_1, \alpha_2\})=\term \bigl(a_2^{-2}a_1\bullet_p a_1a_2^{-1} -a_1a_2^{-2}\bullet_p a_2^{-1}a_1\bigr)= |1|+|-1|=2.$ Thus by Theorem~\ref{maintheorem} any two loops in the free homotopy classes $\alpha_1$ and $\alpha_2$ have at least $2$ intersection points and moreover there are two such loops that do have exactly two intersection points.

{\it Let us now show the Goldman Lie bracket $[\alpha_1, \alpha_2]=0$ and hence it does not allow one to 
compute $\#(\alpha_1, \alpha_2).$\/}

\begin{equation}
\begin{split}
[\alpha_1, \alpha_2]=\langle a_2^{-2}a_1a_1a_2^{-1}\rangle  - \langle a_1a_2^{-2}a_2^{-1}a_1\rangle -
\langle a_1a_2^{-2}a_1a_2^{-1}\rangle + \langle a_2^{-1}a_1a_2^{-1}a_1a_2^{-1}\rangle\\=
\langle (a_1a_2^{-1})(a_2^{-2}a_1^2a_2^{-1})(a_1a_2^{-1})^{-1}\rangle  - \langle a_1a_2^{-3}a_1\rangle \\-
\langle a_1a_2^{-2}a_1a_2^{-1}\rangle + \langle (a_1a_2^{-1})(a_2^{-1}a_1a_2^{-1}a_1a_2^{-1})(a_1a_2^{-1})^{-1}\rangle=0,
\end{split}
\end{equation}
since the first term cancels with the second one and the third term cancels with the fourth one.

{\bf Acknowledgment.} This  work  was  partially  supported  by  a  grant  from  the  Simons  Foundation ($\#235674$ to Vladimir Chernov).

The authors are very thankful to the anonymous referee for the suggestion to include Theorem~\ref{mainthmcor} into the paper.

\end{document}